\documentclass[final,authoryear,5p,times,twocolumn]{elsarticle}

\usepackage{graphicx}
\usepackage{amssymb,amsmath}
\usepackage{multirow}
\usepackage{txfonts}
\usepackage{bm}
\usepackage{lineno}

\newcommand{\ve}[1]{\boldsymbol{#1}}

\begin{document}

\begin{frontmatter}

\title{A blind deconvolution method for ground based telescopes and Fizeau 
interferometers}

\author[Unimore]{M. Prato\corref{cor1}}
\ead{marco.prato@unimore.it}
\author[Unige]{A. La Camera}
\ead{andrea.lacamera@unige.it}
\author[Unife]{S. Bonettini}
\ead{silvia.bonettini@unife.it}
\author[Unimore]{S. Rebegoldi}
\ead{simone.rebegoldi@unimore.it}
\author[Unige]{M. Bertero}
\ead{bertero@unige.it}
\author[Unige]{P. Boccacci}
\ead{patrizia.boccacci@unige.it}
\address[Unimore]{Dipartimento di Scienze Fisiche, Informatiche e Matematiche, Universit\`a di Modena e Reggio Emilia, Via Campi 213/b, 41125 Modena, Italy}
\address[Unige]{Dipartimento di Informatica, Bioingegneria, Robotica e Ingegneria dei Sistemi, Universit\`a di Genova, Via Dodecaneso 35, 16145 Genova, Italy}
\address[Unife]{Dipartimento di Matematica e Informatica, Universit\`a di Ferrara, Via Saragat 1, 44122 Ferrara, Italy}
\cortext[cor1]{Corresponding author. Tel.: +39 059 2055590; fax: +39 059 370513.}

\begin{abstract}
In the case of ground-based telescopes equipped with adaptive optics systems, 
the point spread function (PSF) is only poorly known or completely unknown. 
Moreover, an accurate modeling of the PSF is in general not available. 
Therefore in several imaging situations the so-called blind deconvolution 
methods, aiming at estimating both the scientific target and the PSF from the 
detected image, can be useful. A blind deconvolution problem is severely 
ill-posed and, in order to reduce the extremely large number of possible 
solutions, it is necessary to introduce sensible constraints on both the 
scientific target and the PSF.

In a previous paper we proposed a sound mathematical approach based on a
suitable inexact alternating minimization strategy for minimizing the 
generalized Kullback-Leibler divergence, assuring global convergence.
In the framework of this method we showed that an important constraint on
the PSF is the upper bound which can be derived from the knowledge of its
Strehl ratio. The efficacy of the approach was demonstrated by means of
numerical simulations.

In this paper, besides improving the previous approach by the use of a further
constraint on the unknown scientific target, we extend it to the case of
multiple images of the same target obtained with different PSFs. The main
application we have in mind is to Fizeau interferometry. As it is known this
is a special feature of the Large Binocular Telescope (LBT). Of the two 
expected interferometers for LBT, one, LINC-NIRVANA, is forthcoming while the other, 
LBTI, is already operating and has provided the first Fizeau images, demonstrating the 
possibility of reaching the resolution of a 22.8 m telescope. 
Therefore the extension of our blind method to this imaging modality seems to be timely. 

The method is applied to realistic simulations of imaging both by single 
mirrors and Fizeau interferometers. Successes and failures of the method in 
the imaging of stellar fields are demonstrated in simple cases. These 
preliminary results look promising at least in specific situations. The IDL 
code of the proposed method is available on request and will be included in 
the forthcoming version of the  Software Package AIRY (v.6.1).
\end{abstract}

\begin{keyword}
blind deconvolution \sep numerical optimization \sep image reconstruction 
\sep Fizeau interferometers
\end{keyword}

\end{frontmatter}

\section{Introduction}

In a previous paper \citep{prato2013} we proposed a blind deconvolution method
for ground-based telescopes equipped with an adaptive optics (AO) system. 
Assuming that the image and the corresponding background are known, then the 
features of the method are the following:
\begin{itemize}
\item formulation of the problem as a constrained minimization of the data
fidelity function in the case of Poisson noise (photon counting), namely a
generalized Kullback-Leibler divergence depending on the unknown
astronomical target (in the following called the {\it object}) and on the
unknown point spread function (PSF);
\item non-negativity of the object;
\item non-negativity of the PSF and normalization to unit volume;
\item upper bound on the PSF values derived from the knowledge of the Strehl 
ratio (SR), as suggested by \citet{desidera2009};
\item use of the inexact alternating minimization method proposed by
\citet{bonettini2011} and based on the {\it scaled gradient projection} (SGP) 
method \citep{bonettini2009,prato2012}.
\end{itemize}
The method is iterative and at each outer iteration the object and
the PSF are updated by means of given (but arbitrary) numbers of inner
iterations of SGP. We remark that when SGP is applied to the object only 
projection on the non-negative orthant is required while, in the case of the 
PSF, the projection is performed on the convex and closed set defined by the 
box and equality constraints. 

Our proposed method is similar to blind methods based on the Richardson-Lucy 
(RL) algorithm \citep{richardson,lucy} or its accelerated versions, such as 
those proposed in \citet{holmes,tsumuraya,fish,biggs1998} for the case of 
single image and in \citet{desidera2006} for the case of multiple images, 
applicable to Fizeau interferometry. But the advantage of using SGP in place 
of RL is double: first, as proved in \citet{bonettini2011}, global convergence 
of the iteration holds true, i.e. any limit point of the sequence is a 
stationary point of the constrained KL divergence; secondly it is possible to 
introduce box and equality constraints on the different blocks of variables 
(object and PSF).

Therefore the novelty of the method is that it is based on a sound mathematical 
approach and allows, in an easy way, the introduction of important constraints
on the PSF such as the Strehl constraint, preventing the appearance of trivial 
solutions such as a delta function.

The method, in its present form, does not consider the use of regularization 
for the object or the PSF, in addition to that provided by the constraints 
mentioned above. We remark that it is quite easy to introduce as an additional 
constraint on the object the value of its flux (namely, its $\ell_1$ norm), 
as derived from the detected images and the knowledge of the background. This 
constraint, enforcing the sparsity of the object, is considered in this paper.
If we select a large number of inner iterations on the object variables, the 
method is suitable for the reconstruction of star systems, as already remarked 
in \citet{prato2013}; this result is confirmed in this paper by means of 
simulations more realistic than those used in that paper.

However the main contribution of this paper is the extension of the method to 
the case of Fizeau interferometry. As it is known this is a specific feature 
of the Large Binocular Telescope which consists of two 8.4 m mirrors situated 
on a common mount with a center to center distance of 14.4 m. Indeed, this 
structure is suitable for Fizeau interferometry which should provide images 
with the resolution of a 22.8 m telescope in the direction of the baseline 
joining the center of the two mirrors and that of a 8.4 m telescope in the
orthogonal direction. Different images of the same target corresponding to 
different orientations of the baseline can be combined by suitable 
deconvolution methods to provide a unique reconstructed image with the highest 
resolution in all directions \citep{bertero2011}.

Two interferometers are planned for LBT: the forthcoming LINC-NIRVANA 
\citep{herbst2003}, in advanced realization stage by a German-Italian 
consortium leaded by MPIA, Heidelberg, and the NASA funded LBTI 
\citep{wilson,bailey} already operating on Mount Graham. Indeed, images of 
the Jupiter moon Io were obtained with LBTI/LMIRcam during UT 2013 December 
24, showing that the resolution of a 22.8 m is reachable \citep{leisen} and 
thus proving that LBT is the first in a class of extremely large 
telescopes (ELT).

All the methods developed for Fizeau interferometry are also applicable
to other situations where multiple images of the same target, corresponding to 
different PSFs, are available such as the co-adding problem in Astronomy 
\citep{lucy92} or the multiple image method used in STED microscopy for 
improving the signal-to-noise ratio \citep{castello}.

For this reason we present in Sect. 2 our blind method without a 
specific reference to Fizeau interferometry but just as a method for
multiple image deconvolution in the case of Poisson noise, including, as it 
is obvious, the case of a single image as a particular case. In the same 
section we discuss the intrinsic limitations of our constrained blind 
deconvolution, a discussion which is possible in our sound mathematical 
framework.

In the simulations intended to validate the method we focus on LBT which 
is equipped with a very innovative AO system, the so-called First Light AO 
(FLAO) system \citep{esposito2010}, providing SR values up to 0.9 in K-band. 
Therefore for single image simulations we use models of the PSF of such a 
system. On the other hand, for multiple image deconvolution we consider 
images generated by means of PSFs computed for the interferometer 
LINC-NIRVANA (LN) in K band. Since the camera of LINC-NIRVANA has a pixel 
size of $\sim$ 5 mas while the LMIRcam of LBTI has a pixel size of 10.7 mas, 
the shape of the PSFs of LN in K band (2.2$\mu m$) is similar to the shape 
of the PSFs of LBTI in M band (4.8$\mu m$). Obviously the properties of 
the images may be very different. Details on image modeling and simulation 
are given in Sect. 3. 

Finally in Sect. 4 we discuss our numerical results in the case of binary 
systems and ``open cluster'' models. Conclusions are sketched in Sect. 5. 

\section{Method}

We assume that $p$ different images of the scientific object, with $p$ 
different PSFs are available. The case of a single aperture telescope 
obviously corresponds to $p=1$ if only one image has been acquired; if 
different observations have been performed at different times, hence with
PSFs corresponding to different AO corrections, then the approach can be 
used for the co-adding of these images.

Let ${\ve f}$ be the unknown astronomical object and let ${\ve K}_j$ be the 
unknown PSFs (each one normalized to unit volume) corresponding to the detected 
images ${\ve g}_j$ for $j=1,...,p$ (we assume for simplicity a space-invariant 
model), then we define as $A_j{\ve f}= {\ve K}_j*{\ve f}$ the corresponding 
imaging matrices. Moreover we denote as ${\ve b}_j$ the expected value of the 
background emission in image ${\ve g}_j$ and we assume that it is known, so 
that the expected value of ${\ve g}_j$ is given by $A_j{\ve f}+{\ve b}_j$. 

Since it is quite natural to assume that the $p$ images are statistically 
independent, the likelihood of the problem is the product of the 
likelihoods of the different images. We assume that they are perturbed by
Poisson noise. Then, by taking the negative logarithm of the likelihood we 
obtain the following data-fidelity function which is the sum of $p$ 
Kullback-Leibler (KL) generalized divergences, also known as Csisz{\'a}r 
I-divergences \citep{csiszar}, one for each image, i.e.
\begin{eqnarray}
\label{multipleKL}
& & J_0({\ve f},{\ve K}_1,...,{\ve K}_p ;{\ve g},{\ve b})= \\ 
\nonumber
& & \sum_{j=1}^{p}\sum_{{\bf m} \in S}\Big\{{\ve g}_j({\bf m}) {\rm ln} 
\frac{{\ve g}_j({\bf m})}{(A_j \ve{f})({\bf m})+{\ve b}_j({\bf m})}~+ \\ 
\nonumber
& & (A_j \ve{f})({\bf m})+{\ve b}_j({\bf m})-{\ve g}_j({\bf m}) \Big\}~~,
\end{eqnarray}
where $S$ is the set of the values of the multi-index ${\bf m}$ 
characterizing the pixels of the image array, and
$({\ve g},{\ve b}) = \{({\ve g}_j,{\ve b}_j)\}_{j=1}^p$. The problem 
of image deconvolution (without regularization) consists in the minimization 
of this function with respect to ${\ve f}$ for given PSFs, images and 
backgrounds. The minimization can be obtained by means of RL in the 
single image case, by means of OS-EM method \citep{hudson,bertero2000} in the
multiple image case, or by means of SGP method \citep{bonettini2009,prato2012}
in both cases. As shown in \citet{prato2012} SGP is more efficient than 
OS-EM if the number $p$ of images is not too large.

\subsection{Blind deconvolution: problem formulation}

If images and backgrounds are given, the ML approach to blind deconvolution 
can be formulated as the minimization of the function in Eq. (\ref{multipleKL}) 
with respect to $p+1$ blocks of unknown variables, namely the object
${\ve f}$ and the $p$ PSFs ${\ve K}_j$, $j=1,...,p$. As it is known, 
this function is convex with respect to each block of variables for fixed
values of the others, but is not convex with respect to the full set of
variables \citep{prato2013}. Therefore blind deconvolution is a difficult
problem of nonconvex optimization. Moreover, this problem is highly 
ill-posed and allows uninteresting solutions. For instance, a global minimum 
can be achieved by choosing ${\ve f}={\ve{g - b}}$ and 
${\ve K}_j = \bm{\delta}$, $j=1,\ldots,p$, where $\bm{\delta}$ is the
Dirac delta array. This trivial minimizer can be avoided by
the introduction of suitable regularization terms and constraints.

Since we mainly consider the case of stellar fields or, in other words, of
sparse objects, by taking into account the sparsity property of the
minimizers of the KL divergence in the case of image deconvolution
\citep{bertero2009}, we do not introduce an object-dependent regularization
term in the objective function. However, besides non-negativity of the
object we also introduce a constraint on its flux; more precisely we require
that the object flux coincides with the average flux of the $p$ detected
images (after background subtraction), which is given by
\begin{equation}
c = \frac{1}{p}\sum_{j=1}^p\sum_{{\bf m} \in S}\{{\ve g}_j({\bf m})-{\ve b}_j
({\bf m})\}~~.
\end{equation}
We remark that this constraint is further enforcing sparsity; in the case of 
deconvolution and zero value of the backgrounds, it is automatically 
satisfied by the minimizers of the KL divergence.
 
As concerns the PSFs, as shown in \citet{desidera2009} and 
\citet{prato2013}, an important constraint is the upper bound derived from
the knowledge of the SR characterizing the AO correction of the
atmospheric blur during the observation. Moreover, non-negativity and
normalization provide additional constraints. In conclusion, the nonconvex
optimization problem we are considering can be formulated as follows 
\begin{eqnarray}
\label{prob}
& &{\rm min}~~J_0({\ve f},{\ve K}_1,...,{\ve K}_p;{\ve g},{\ve b}) \\ 
\nonumber
& &{\rm s.t.}~~~{\ve f}\geq 0~,~\sum_{{\bf n}\in S}{\ve f}({\bf n})=c~; \\
\nonumber
& & ~~~~~~~0\leq {\ve K}_j\leq s_j~,\sum_{{\bf n}\in S} {\ve K}_j({\bf n}) = 
1~;~j=1,...,p~~,
\end{eqnarray}
where $s_j$ is the upper bound on the PSF ${\ve K}_j$ derived from the 
knowledge of the SR characterizing the acquisition of ${\ve g}_j$. In
conclusion the data of the problem are $({\ve g}, {\ve b})$ and
$s = \{s_j \}_{j=1}^p$.

Another important constraint can be provided by the requirement of 
band-limiting of the PSFs, which is used for instance in \citet{desidera2006}.
Indeed the band of the PSF, i.e the set in Fourier space where the
Fourier transform of the PSF is not zero, is known and consists, in general,
of a disc in the case of a single aperture telescope and a union of three discs
in the case of a Fizeau interferometer with LBT, the central one being the 
band of the mirrors of LBT and the side-ones the replicas due to interferometry 
\citep[see, for instance,][]{bertero2011}. This constraint on the Fourier 
transform of the PSF, together with normalization and upper and lower bounds, 
defines a convex set. Unfortunately, since we use projection methods, the 
projection of an array on this set is not easily computable, even if methods for 
computing the projection on the intersection of convex sets are available. 
The difficulty is that these methods are not efficient and therefore can lead
to an excessive computational cost since the projection should be computed
several times in the used iterative methods. However, from the numerical 
experiments described in \citet{prato2013} we deduce that a suitable 
initialization of our algorithm, based on a PSF satisfying the band constraint,
may lead to a set of reconstructed PSFs whose bands are very close to the 
desired ones.

\subsection{Blind deconvolution: alternating minimization}

Although the previous formulation of blind deconvolution requires the 
minimization of a nonconvex objective function, the constraints have a nice 
separable structure, since they involve separately the blocks of variables, 
defining a feasible convex set for each of them. In addition, the function is 
convex with respect to the different blocks of variables. In these settings, 
the solution of problem \eqref{prob} can be sought by means of an 
{\it alternating minimization} (AM) strategy.

The basic idea of AM is the cyclic minimization for the constrained problem
with respect to one block of variables, updating its value for the next
minimization step. This kind of approach is known in the literature also as 
{\em nonlinear Gauss-Seidel} or {\em block coordinate descent} method and its
theoretical properties have been deeply studied in the last decades 
\citep{bertsekas,Grippo1999,Grippo2000,Luo1992}.

In our case, each iteration of the AM method consists in solving the following 
$p+1$ constrained minimization problems of convex functions
\begin{eqnarray}
\label{am}
& &{\ve f}^{(k+1)} = {\rm arg~min}_{{\ve f} \in \Omega}J_0({\ve f},{\ve K}_1^{(k)},...,{\ve K}_p^{(k)};{\ve g},{\ve b}) \\ \nonumber
& &{\ve K}_1^{(k+1)}={\rm arg~min}_{{\ve K} \in \Omega_1}J_0({\ve f}^{(k+1)},{\ve K},...,{\ve K}_p^{(k)};{\ve g},{\ve b}) \\ \nonumber
& &\quad \vdots \\\nonumber
& &{\ve K}_p^{(k+1)}={\rm arg~min}_{{\ve K} \in \Omega_p}J_0({\ve f}^{(k+1)},{\ve K}_1^{(k+1)},...,{\ve K};{\ve g},{\ve b})~~, 
\end{eqnarray}
where $k$ is an index running on the AM iterations, $\Omega$ is the set of 
the constraints on the object ${\ve f}$ and $\Omega_j$ is the set of the 
constraints on the PSF ${\ve K}_j$. The limit points of the sequence 
generated by this iteration scheme are also stationary points for the constrained problem 
if each partial problem is the minimization of a strictly convex 
function \citep{bertsekas}. This condition is not only sufficient but also 
necessary when more than two blocks of variables are involved. Indeed, 
\citet{powell} showed a counterexample where the strict convexity is not 
satisfied, three blocks of variables are involved and the AM method fails to 
locate stationary points.

Even when the hypothesis of strict convexity on each block of variables holds 
true, the convergence of the AM scheme can be proved only if each partial 
minimization problem is solved {\it exactly}, which is often impractical or 
too costly (in general iterative methods are used). Quite surprisingly the 
convergence result is obtained without the assumption of strict convexity if 
SGP is used for solving inexactly the partial minimization problems. This 
important result is proved in \citet{bonettini2011} and is basic for the 
proposed approach to blind deconvolution.

In conclusion, at each line of Eq. (\ref{am}) the minimization is replaced
by a given number of SGP iterations. These will be called {\it inner}
iterations while the iterations of the AM scheme will be called {\it outer}
iterations. Therefore the sequence generated by the method depends both
on the initialization of the first outer iteration (in the subsequent outer 
iterations the inner iterations are initialized with the results derived 
from the previous one) and on the given numbers of inner iterations for 
each block of variables. In the present application it seems quite natural 
to choose the same number of inner iterations for all the PSF's blocks.

\subsection{SGP algorithm}

Since it is basic for the solution of the partial minimization problems
it may be useful to briefly recall the main points of the SGP
algorithm \citep{bonettini2009} even if its application to astronomical 
imaging has already been described elsewhere \citep{prato2012,Bonettini2010,
Bonettini2014}. 
To this purpose we remark that each minimization problem in the iterative 
scheme of Eq. (\ref{am}) has the following structure
\begin{equation}
\min_{\ve h \in \Omega} J_0({\ve h})~~,
\label{minpr}
\end{equation}
where, for simplicity, we omitted the dependence on the other variables and
$\Omega$ is the closed and convex set defined by the constraints. The main 
difference with respect to \citet{prato2012} is that $\Omega$ is a subset of 
the non-negative orthant defined by a suitable equality constraint. Therefore 
the projection on this set is more complex than that on the non-negative orthant.

The main step of SGP is the computation of the $k$-th feasible descent direction 
(where $k$ is an index running on the inner iterations of a given AM iteration)
\begin{equation*}
{\ve d^{(k)}} = P_{\Omega,D_k^{-1}}({\ve h^{(k)}}-\alpha_kD_k
\nabla J_0({\ve h}^{(k)})) - {\ve h}^{(k)}
\end{equation*}
by performing the following steps:
\begin{itemize}
\item[a)] The direction provided by the negative gradient $-\nabla J_0({\ve h}^{(k)})$
is modified by a diagonal scaling matrix $D_k$ with positive entries, which in
all the subproblems of one AM iterations is given by
\begin{equation}
D_k = \mbox{diag}\left(\min(L_2,\max(L_1,{\ve h}^{(k)})\right)~~,
\end{equation}
$(L_1,L_2)$ being given constants estimated from the extreme values of the image.
\item[b)] A point on the scaled gradient direction is selected by choosing a multiplicative
factor $\alpha_k$ by means of an alternation of the generalized Barzilai-Borwein 
(BB) rules \citep{barzilai,bonettini2009}
\begin{eqnarray}
\label{BB1-2}
& & \alpha_k^{{(BB1)}} = \frac{({\ve s}^{(k-1)})^T D_k^{-1} D_k^{-1}
{{\ve s}^{(k-1)}}}{({\ve s}^{(k-1)})^T D_k^{-1} {\ve z}^{(k-1)}}~~, \\ \nonumber
& & \alpha_k^{{(BB2)}} = \frac{({\ve s}^{(k-1)})^T D_k {\ve z}^{(k-1)} }
{({\ve z}^{(k-1)})^T D_k D_k {\ve z}^{(k-1)}}\,,
\end{eqnarray}
where ${\ve s}^{(k-1)}={\ve h}^{(k)}- {\ve h}^{(k-1)}$ and
${\ve z}^{(k-1)}=\nabla J_0({\ve h}^{(k)})-\nabla J_0({\ve h}^{(k-1)})$, and
a suitable introduction of upper and lower bounds.
\item[c)] The resulting point is brought back in the feasible set $\Omega$ by 
means of the 
projection $P_{\Omega,D_k^{-1}}$ associated to the norm induced by $D_k^{-1}$, i.e.
\begin{equation}
P_{\Omega,D_k^{-1}}({\ve h})= {\rm arg~min}_{{\ve y}\in\Omega} 
({\ve h-\ve y})^TD_k^{-1}({\ve h-\ve y}).
\label{proj}
\end{equation}
Since the feasible sets of both the object and the PSFs involve a given number 
of inequalities plus 
an equality constraint, in all cases we used a secant-based routine developed by \citet{Dai2006}, 
which is able to compute the projection with a computational cost growing linearly in time with 
respect to the image size \citep[see also][]{prato2013}. 
\end{itemize}

\subsection{Discussion}

In the previous approach, the blind deconvolution problem is formulated as the
constrained minimization of a nonconvex function which depends on an
extremely large number of variables, about $10^6$ in the numerical
experiments described in this paper. Since the constraints used in our
approach imply that the sequences of objects and PSFs generated by the inexact
AM method are bounded, it follows that these sequences have limit points.
We can add that, even if it is difficult to provide a theoretical evidence of
the existence of a unique limit point, in all our numerical experiments the 
sequences produced by the inexact AM method have a convergent behaviour.
However, according to the general convergence result proved in 
\citet{bonettini2011}, we can only state that the limit points are 
{\it stationary points} of the function, hence not necessarily minimizers. 

As far as we know, there is no practical way for establishing if these points 
are minimizers or not. In fact, it should be necessary to manage the Hessian 
of the function in these points but this is an absolutely intractable matrix 
even if one can write it explicitly \citep{prato2013}.
Since one can use different initializations of the iterative procedure and
different numbers of inner iterations and these different choices can 
produce different results, in a practical application we do not see an 
approach better than that of doing different attempts and look for that 
providing the most sensible solution.

An additional difficulty is that it may happen, as we show by some
numerical simulations, that a sensible solution corresponds to a value
of the objective function which is greater than the value of the same
function corresponding to a solution which is clearly unphysical. It is
obvious that these situations should not be surprising because the problem
of blind deconvolution is nonconvex and therefore the objective function
can have several local minimizers as well as stationary points. Since the
objective function has a simple structure it should be important to
characterize the sets of these points, but an approach to this problem
presently is not available, as far as we know.

The advantage of our method is that it is mathematically sound, it provides
sequences with limit points, very frequently with a unique limit point and
therefore, if the user is conscious of the difficulties of the problem,
he can attempt to use this method for obtaining different solutions in
practical applications and select that looking as the most appealing. In 
the next section we attempt to provide a few hints for helping the user in
the choice of the parameters of the method and, in particular, of its
initialization.

\section{Image simulation}

We model the images according to the model proposed in \citet{snyder1994} 
for images acquired with a CCD camera, i.e. each pixel is affected by 
background (due to sky emission, dark current, etc.), photon counting noise
(described by a Poisson distribution) and additive read-out noise (RON) 
described by a Gaussian distribution. 

If the RON variance is $\sigma^2$, in the deconvolution process it can be 
approximated by a Poisson distribution with parameter $\sigma^2$ if $\sigma^2$
is added both to the detected images and the corresponding backgrounds 
\citep{snyder1995}. Therefore all the pixel values of the detected images 
can be viewed as realizations of suitable Poisson random variables if in 
Eq. (\ref{multipleKL}) we intend that ${\ve g}_j,{\ve b}_j$ have been 
modified according to this approach. Therefore, in our numerical simulations 
we perturb the images with Poisson and additive Gaussian noise but in the 
deconvolution algorithms we use the images and backgrounds modified as above.

All the images and the PSFs considered in our numerical experiments are sized 
$256 \times 256$ pixels in the single image case, with a pixel size of 15 mas, 
and $512 \times 512$ pixels in the multiple image case, with a pixel size of
5 mas. Moreover all images, except one indicated in Sect. 4.1.1, are obtained 
by adding 10 frames in order to avoid saturation of the detector, as we 
discuss in the following, so that the variance of the RON will be $10~\sigma^2$.

\subsection{Single image simulation}

In this case we use two PSFs in K-band with SR = 0.81 and 0.62 respectively, 
modeling the optics of a single mirror of LBT, with diameter 8.4 m, and the 
effect of the adaptive optics system FLAO using the power spectrum of the 
wavefront residual of the AO correction as measured at the telescope 
\citep{esposito2012}. To the noise-free image, obtained by convolving the 
object with one of these PSFs, a background in K-band is added and the result 
is corrupted with Poisson and additive Gaussian noise. In order to avoid 
saturation of the detector (a maximum number of $5 \times 10^4$ photons per 
pixel is assumed in a single frame) the image is obtained by co-adding $n$ 
frames. More precisely, in the case of a stellar system the procedure for 
image generation is the following.
\begin{itemize}
\item We establish the coordinates of the stars and we fix their 
magnitudes in K-band.
\item We compute the integration time which does not produce saturation of the 
detector by taking into account the collection area of the telescope, 
the overall efficiency of the acquisition system (assumed equal to 30\%), 
and the flux of the brightest star multiplied by the peak value 
of the PSF. This is the integration time of a single frame and is used for 
computing the number of frames $n$ required for obtaining an acceptable SNR 
for all the stars of the system.
\item We generate noise-free images by shifting, with sub-pixel precision, 
the PSF to the positions of the stars and adding these shifted PSFs, each one 
weighted with a weight corresponding to the magnitude and the total observation 
time.
\item These images are perturbed by adding a background in K-band, 
corresponding to about 13.5 mag arcsec$^{-2}$, and by corrupting the
results with Poisson and additive Gaussian noise (RON); the variance of the
RON is $n \sigma^2$, thus corresponding to the RON of $n$ frames; we take
$\sigma = 10~e^-/px$.
\end{itemize}

\subsection{Multiple image simulation}

As concerns the simulation of LN images, we recall that the instrument combines 
in a Fizeau mode the beams coming from the two mirrors of LBT whose 
center-to-center distance is about 14.4 m. Therefore the maximum baseline 
available is 22.8 m and the resolution achievable by a single LN image is 
that of a 22.8 m telescope in the direction of the baseline and that of a 8.4 m 
telescope in the orthogonal direction. For a given orientation the PSF of LN
looks as that of a 8.4 m telescope, modulated by the interference fringes,
orthogonal to the direction of the baseline. In order to get a more uniform 
resolution one must acquire and combine different images with different 
orientations of the baseline.

It is important to remark that the orientation of the fringes does not depend 
on the orientation of the baseline because the camera is rotating with the 
baseline and therefore the fringes have always the same direction (for 
instance the vertical one) in the image array. In other words two images of 
the same scientific object with two different orientations of the baseline 
correspond to two rotated versions of that object. 
This specific feature implies that one should introduce rotation matrices in 
the formulation of the problem. However we verified that the computation 
of hundreds or thousands of rotations in hundreds or thousands of inner 
iterations introduces large computational errors. Therefore we considered 
the approach which consists in derotating the images in such a way that they 
correspond to aligned versions of the object $\ve f$. The price to be payed 
is that the derotation of discrete images modifies their statistical properties. 
In order to estimate this effect we considered the rotation of a constant 
array perturbed by Poisson noise. We found the following results:
\begin{itemize}
\item before rotation the histogram of the array is a Gaussian with the same
mean and variance; after a rotation based on spline interpolation the 
histogram is still a Gaussian with the correct mean but a smaller variance;
\item the support of the autocorrelation of the rotated image is a 
$3 \times 3$ square;
\item if we use a different rotation approach which consists in attributing
the value of a pixel before rotation to the pixel with maximum overlapping
after rotation (nearest neighbor approximation), the statistics 
is preserved but the quality of the image is degraded.
\end{itemize}
As a consequence of this analysis we decided to use in the approach derotated
images. 

The procedure adopted in our numerical experiments is similar to that used in 
the case of a single image. We consider two sets of PSFs in K-band with SR 
respectively 0.77 and 0.46, corresponding to orientation angles of the baseline 
indicated as $0^\circ,~60^\circ$ and $120^\circ$, all with vertical fringes (for
simplicity we take the same SR for the three orientations). The 
first PSF of each set has been generated by means of the software package LOST 
\citep{arcidiacono2004}, the second by reflecting the first one with respect 
to the central line and the third by taking the arithmetic mean of the first 
two. In this way the three PSF of each set have exactly the same SR. 
Then the generation of the corresponding LN images is similar to that of
the single image case by modifying the first item as follows. 
\begin{itemize}
\item We establish the coordinates of the stars corresponding to the
observation at $0^\circ$ and we compute, with sub-pixel precision, their 
coordinates if the system is rotated by $60^\circ$ and $120^\circ$ respectively.
\end{itemize}
The rest of the procedure is unchanged and applied to the three images but at
the end we must add the following item.
\begin{itemize}
\item The images corresponding to $60^\circ$ and $120^\circ$ are
derotated in order to align the object in the three images and three
arrays containing the object are extracted from the full images.    
\end{itemize}
The derotated images are used in the definition of the objective function 
and in the blind algorithm, which therefore will produce derotated PSFs. 

\section{Numerical results}

In order to evaluate the quality of the reconstructions obtained with our 
blind method we need some figures of merit.

As concerns the reconstruction of a binary we consider the relative absolute
error on the magnitudes of both stars while in the case of a stellar system
we consider a magnitude average relative error (MARE) defined by
\begin{equation}
\label{mare}
MARE=\frac{1}{q}\sum_{i=1}^q\frac{|m_i-\tilde m_i|}{\tilde m_i}~~,
\end{equation}
where $q$ is the number of stars and $m_i$, $\tilde m_i$ are respectively 
the reconstructed and the true magnitudes.

As concerns PSF reconstruction, in the case of single image we consider the 
root-mean-square error with respect to the true one, defined as usual in terms 
of the $\ell_2$ norm of their difference. In the
case of LN images generated according to the previous procedure, since
the blind algorithm produces a set of three PSFs, two of them being derotated 
with respect to the ones used for generating the images, for comparison 
we must derotate the original ones. If we denote as $\tilde{K}_j$ the
derotated original PSF, then we measure the quality of the reconstruction by 
means of the root-mean-square error (RMSE)
\begin{equation}
\rho_j=\frac{\|K_j - \tilde{K}_j\|}{\|\tilde{K}_j\|}~~,
\end{equation}
where $K_j$ is the reconstructed PSF and $\|\cdot\|$ denotes the usual 
$\ell_2$-norm. 

\subsection{Binary systems}

We first consider the simple case of binary systems. More precisely we 
consider nine cases by varying both separation and magnitude of the stars. 
By keeping fixed the magnitude of the primary, i.e. $m_1 = 15$, we take for 
the magnitude of the secondary $m_2 = 15,~16$ and 17. Moreover for each 
choice we consider three possible angular separations: $d = 60$, $120$ and
$240$ mas in the single image case and $d = 20$, $40$ and $80$ mas in the LN 
case. In both cases the first separation corresponds to the resolution limit 
of the instrument while the last is four times larger.
In all cases, as described in the previous section, 
we compute the integration time of a frame in such a way that the number 
of counts in the image pixel corresponding to the position of the primary
does not exceed $5 \times 10^4$. As stated in the previous section, we 
consider 10 frames per image, both in the single and in the multiple image 
case, so that the peak value of the photons is about $5 \times 10^5$ for all 
images. Since in the case of LN we have three images, in this case the SNR is 
higher than in the single image case.

In Fig. 1 we show the images of the binaries with 
$m_1 = m_2 = 15$ and different angular separations; in the
first row those of the single image case and in the second row those
of the multiple image case corresponding to the $0^\circ$ baseline, all
obtained with the PSF with the highest SR. The difficulty in reconstructing 
the binary with separation $d = 20$ mas is obvious. 

\begin{figure*}
\centering
\includegraphics[width=0.6\textwidth]{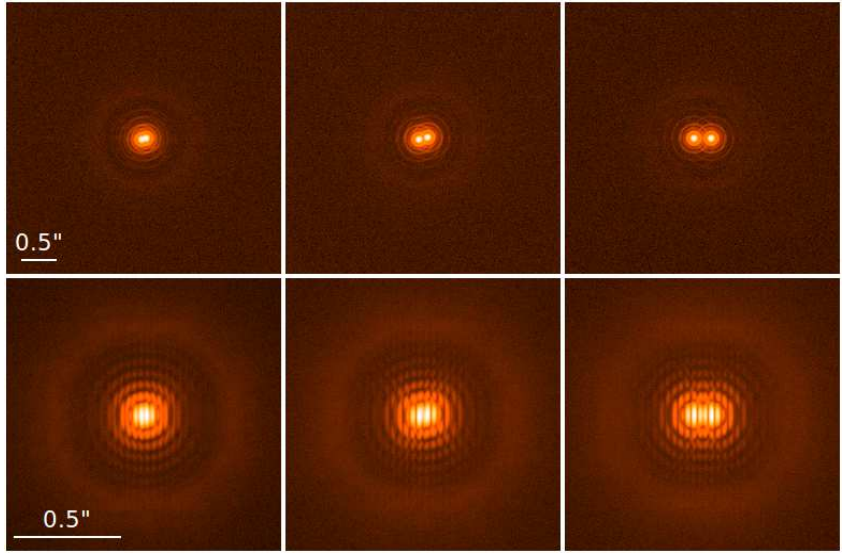} 
\caption{Examples of input images of binaries with magnitudes $m_1 = m_2 = 15$.
In the first row those of the single image case, corresponding to PSF with 
SR = 0.81: from left to right, angular separation of 60, 120 and 240 mas. 
In the second row those of the multiple image case, corresponding to the
PSF with SR = 0.77: from left to right, angular separation of 20, 40 and 80 
mas. These images correspond to the first orientation of the baseline and only
the central part of the images 256 $\times$ 256 is displayed. In the 
two other orientations the binaries appear rotated by 60 and 120 degrees 
respectively. Images are displayed in log scale. The length corresponding to 
0.5 arcsec is also indicated.}
\label{fig:binaries}
\end{figure*}  

\subsubsection{Single image}

For the convenience of the reader we give the computed integration time
avoiding saturation in a single frame: 40 sec for SR = 0.81 and 52 sec for 
SR = 0.62. As already stated the images are obtained by adding 10 frames. 
These are the input images of the blind algorithm together with the value of 
the background.

\begin{figure*}
\centering
\includegraphics[width=\textwidth]{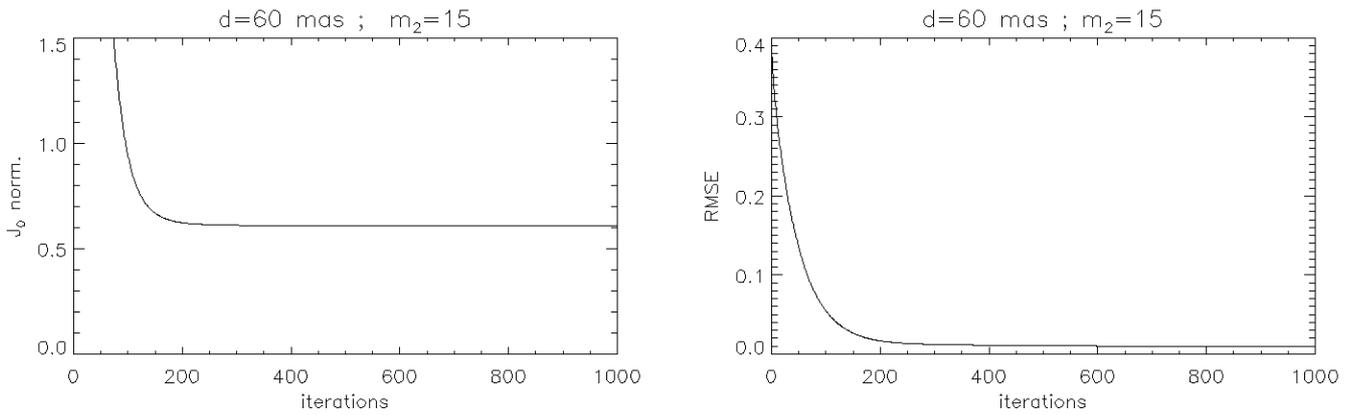}
\caption{Behaviour, as a function of the number of iterations, of the
normalized objective function (left panel) and of the RMSE on the
PSF (right panel). The parameters of the binary are indicated 
in the figure. The plots refer to the PSF with SR = 0.81.}
\label{fig:objective}
\end{figure*}

In a first attempt we use the initialization already used in \citet{prato2013} 
and in other papers, namely a constant array for the object and the 
autocorrelation of the diffraction-limited PSF for the PSF. Indeed, this 
initialization has produced very promising results in our previous paper,
where a much higher SNR was assumed. We use 1000 outer iterations in the 
case SR = 0.81 and 2000 outer iterations in the
case SR = 0.62. Indeed in the case of a lower SR we have a lower quality of
the images and, presumably, a larger number of iterations is required. As
concerns the inner iterations, as in \citet{prato2013} we use 50 SGP
iterations for the object and one SGP iteration for the PSF.

In Table 1 we give the results obtained with the previous choice.
As a first remark, the binaries and the PSFs are reconstructed satisfactorily
in all cases except the closest binaries ($d = 60$ mas) with different
magnitudes. Indeed, the indication 100\% in the column for $\Delta m_2/m_2$ 
means that the method reconstructs only one star, which sometimes is not 
exactly in the position of the primary but slightly shifted in the direction 
of the secondary. Since its magnitude is computed using a 3$\times$3 square 
centered on the true position of the primary, the error on its magnitude is, 
in general, not too large. On the other hand the error on the PSFs is very 
large, as one should expect since the secondary is missed. This point deserves 
further investigation.

In Fig. 2 we show, in a particular case, the behaviour of the 
normalized objective function, defined by $2 J_0/ N^2$ with $J_0$ given in 
Eq. (\ref{multipleKL}) (with $p = 1$), and of the RMSE on the PSF as functions of
the number of iterations. Similar behaviors are obtained in all cases where
a sensible result is obtained. This result suggests that presumably 
convergence is reached after 1000 iterations even if, as previously discussed, 
it is difficult to establish numerically the convergence of a sequence.

\begin{table*}
\caption{Single image case - Binary reconstructions provided by the algorithm 
initialized with the autocorrelation of the diffraction-limited PSF. In the 
first column the value of the SR, in the second the angular separation, in 
the third the magnitude of the secondary, in the fourth and fifth the errors 
on the magnitudes of the two stars. In the subsequent column we give the RMSE 
for the reconstructed PSF. Finally in the last two columns we give the value
of the normalized objective function, defined by $2J_0/N^2$, as computed
at the end of the iterations, and the number of outer iterations.}
\label{binaries_s1}
\centering
\begin{tabular}{|c|c|c|cc|c|c|c|}
\hline
SR & d (mas) & $m_{2}$ & $\Delta m_1/m_1$ & $\Delta m_2/m_2$ & RMSE & $J_0$ norm. & IT \\
\hline
  \multirow{9}{*}{0.81} & \multirow{3}{*}{60} & 15 & 0.05\% & 0.03\% & 0.94\% & 0.6071 & 1000 \\
  & & 16 & 2.37\% & 100\% & 40.41\% &  0.5549 & 1000 \\
  & & 17 & 0.99\% & 100\% & 16.77\% &  0.5964 & 1000 \\
\cline{2-8}
  & \multirow{3}{*}{120} & 15 & $<$0.01\% & $<$0.01\% & 0.76\% & 0.6047 & 1000 \\
  & & 16 & 0.03\% & $<$0.01\% & 1.11\% & 0.6296 & 1000 \\
  & & 17 & 0.03\% & 0.17\% & 1.40\% & 0.6254 & 1000 \\
\cline{2-8}
  & \multirow{3}{*}{240} & 15 & $<$0.01\% & $<$0.01\% & 0.79\% & 0.5999 & 1000 \\
  & & 16 & 0.02\% & 0.03\% & 0.83\% & 0.6273 & 1000 \\
  & & 17 & $<$0.01\% & 0.04\% & 1.17\% & 0.6229 & 1000 \\
\hline
  \multirow{9}{*}{0.62} & \multirow{3}{*}{60} & 15 & 0.18\% & $<$0.01\% & 1.08\% & 0.5338 & 2000 \\
  & & 16 & 2.31\% &  100\% & 34.37\% & 0.4635 & 2000 \\
  & & 17 & 1.05\% &  100\% & 16.87\% & 0.4983 & 2000 \\
\cline{2-8}
  & \multirow{3}{*}{120} & 15 & 0.15\% & 0.14\% & 1.04\% & 0.5261  & 2000 \\
  & & 16 & 0.02\% & 0.01\% & 1.28\% & 0.5419 & 2000 \\
  & & 17 & 0.04\% & 0.25\% & 1.59\% & 0.5329 & 2000 \\
\cline{2-8}
  & \multirow{3}{*}{240} & 15 & 0.04\% & 0.04\% & 1.00\% & 0.5309 & 2000 \\
  & & 16 & $<$0.01\% & 0.06\% & 1.13\% & 0.5537 & 2000 \\
  & & 17 & 0.05\% & 0.36\% & 1.80\% & 0.5361 & 2000 \\
  \hline
\end{tabular}
\end{table*}

A second remark is that, according to statistical properties of Poisson random 
variables, if we compute the value of the normalized objective function by 
inserting in Eq. (\ref{multipleKL}) the noisy and the noise-free images  
we should obtain a value very close to 1 \citep{bertero2010,zanella2009}. This 
is just what we obtain using our simulated images (this result also demonstrates 
the accuracy of the approximation of the RON with a Poisson random variable). 
However the limiting values of the normalized objective function obtained in our 
experiments are definitely smaller than 1, an effect already remarked in our 
previous paper.

Coming back to the problem of the unresolved binaries, we point out that, if 
we deconvolve the images using the PSF used for their generation 
({\it inverse crime}) all the binaries are correctly reconstructed with small 
errors on their magnitudes. Therefore the failure of our experiment may be 
due to a failure of the method or to an inappropriate initialization or to 
inappropriate choices of the internal iterations. 

Several attempts with different numbers of internal iterations did not improve 
the results. Therefore we searched for an initial PSF with a SR value closer 
to the correct one and with the property of being band-limited with the band 
of the LBT mirror. A possible choice is obtained by means of the 
diffraction-limited PSF of LBT, let us say $\tilde K$, by looking for an
initial guess $K^{(0)}$ of the following form
\begin{equation}
K^{(0)}= \frac{1}{1 + \omega~N^2}(\tilde K + \omega)
\end{equation}
which is band-limited and satisfies the normalization condition. The constant 
$\omega$ should be selected in such a way that $K^{(0)}$ has the correct
SR value, i.e. $\max~(K^{(0)})$ = SR $\max~(\tilde K)$. We obtain
\begin{equation}
\big({\rm SR} ~ N^2~\max(\tilde K) - 1\big) \omega = 
({\rm SR} - 1)~\max(\tilde K)
\end{equation}
and, by neglecting 1 with respect to the first term in the l.h.s. of this
equation, we obtain $\omega=(1-{\rm SR})/({\rm SR}~ N^2)$. 

The results obtained 
with this initialization, using again 50 SGP iterations for the object and 
one for the PSF, are reported in Table 2. Since the 
convergence is slower than in the previous case we use 2000 outer iterations 
for SR = 0.81 and 3000 iterations for SR = 0.62.

\begin{table*}
\caption{Single image case - Binary reconstructions provided by the algorithm 
initialized with the diffraction-limited PSF plus a constant selected for 
satisfying the SR constraint (see the text). The structure of the Table is the 
same of Table \ref{binaries_s1}.}
\label{binaries_s2}
\centering
\begin{tabular}{|c|c|c|cc|c|c|c|}
\hline
SR & d (mas) & $m_{2}$ & $\Delta m_1/m_1$ & $\Delta m_2/m_2$ & RMSE & $J_0$ norm. & IT \\
\hline
  \multirow{9}{*}{0.81} & \multirow{3}{*}{60} & 15 & 0.02\% & 0.04\% & 0.82\% & 0.6067 & 2000 \\
  & & 16 & 0.09\% & 0.16\% & 2.05\% & 0.6237 & 2000 \\
  & & 17 & 1.01\% & 100\% & 17.08\% & 0.5956 & 2000 \\
\cline{2-8}
  & \multirow{3}{*}{120} & 15 & $<$0.01\% & $<$0.01\% & 0.77\% & 0.6046 & 2000 \\
  & & 16 & 0.02\% & $<$0.01\% & 1.09\% & 0.6294 & 2000 \\
  & & 17 & 0.02\% & 0.15\% & 1.35\% & 0.6253 & 2000 \\
\cline{2-8}
  & \multirow{3}{*}{240} & 15 & $<$0.01\% & $<$0.01\% & 0.80\% & 0.5989 & 2000 \\
  & & 16 & 0.02\% & 0.02\% & 0.82\% & 0.6271 & 2000 \\
  & & 17 & $<$0.01\% & 0.02\% & 1.12\% & 0.6227 & 2000 \\
  \hline
  \multirow{9}{*}{0.62} & \multirow{3}{*}{60} & 15 & 0.02\% & $<$0.01\% & 1.11\% & 0.5333 & 3000 \\
  & & 16 & 0.12\% & 0.25\% & 2.64\% & 0.5354 & 3000 \\
  & & 17 & 1.05\% & 100\% & 16.87\% & 0.4983 & 3000 \\
\cline{2-8}
  & \multirow{3}{*}{120} & 15 & 0.01\% & 0.01\% & 1.06\% & 0.5258 & 3000 \\
  & & 16 & 0.02\% & $<$0.01\% & 1.26\% & 0.5419 & 3000 \\
  & & 17 & 0.04\% & 0.25\% & 1.58\% & 0.5329 & 3000 \\
\cline{2-8}
  & \multirow{3}{*}{240} & 15 & 0.03\% & 0.03\% & 0.99\% & 0.5304 & 3000 \\
  & & 16 & $<$0.01\% & 0.06\% & 1.12\% & 0.5537 & 3000 \\
  & & 17 & 0.05\% & 0.36\% & 1.80\% & 0.5361 & 3000 \\
  \hline
\end{tabular}
\end{table*}

\begin{figure*}
  \centering
  \includegraphics[width=0.6\textwidth]{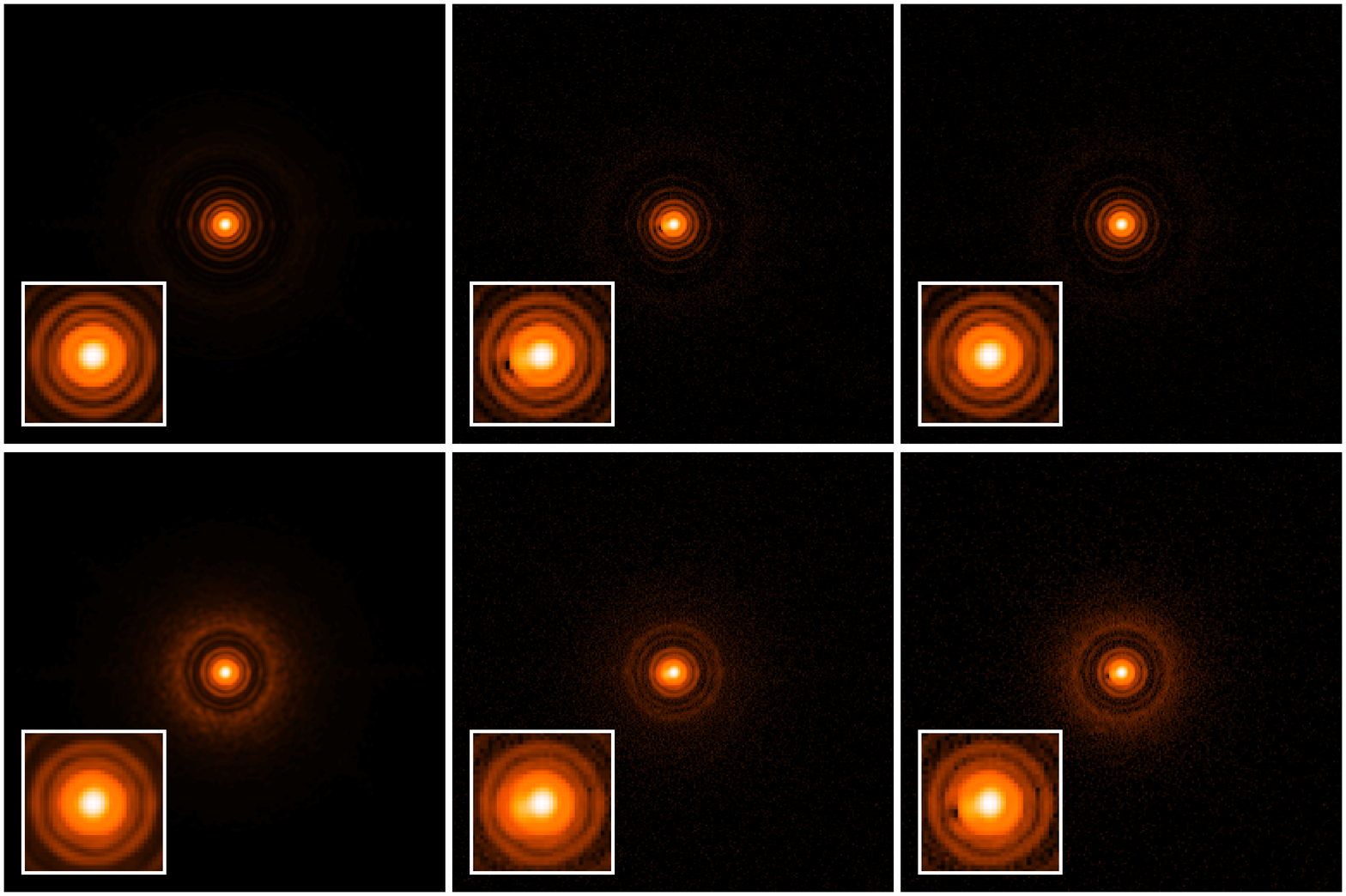} 
  \caption{Single image case - PSF reconstruction in the case of the binary with 
$d = 60$ mas and $m_2=17$. The input image is the sum of 30 frames (see text). 
These PSFs correspond to the minima of the reconstruction error.
First column: the true PSF with SR = 0.81 (top) and SR = 0.62 (bottom). 
Second column: PSF reconstruction provided by the algorithm initialized with 
the autocorrelation of the diffraction-limited PSF. Last column: PSF 
reconstruction provided by the algorithm initialized with the 
diffraction-limited PSF plus a constant. In each panel we also show a zoom 
of the core of the PSF which makes evident artifacts due to the secondary. 
All  images are displayed in log scale.}
\label{fig:psf_s}
\end{figure*}

By comparing the results reported in the two tables we remark that the
two different initializations provide  very similar results in all cases
where they succeed or they fail; in other words they provide sequences of
iterations which presumably converge, even if with a different rate, to
the same point, which is a stationary point of the objective function. 
Obviously we believe that it is also a minimizer. In the case of separation 
60 mas and $m_2 = 16$ the algorithm, equipped with the new initialization, 
is able to reconstruct the binary and the PSF with a satisfactory accuracy for
both values of SR. We remark that the value of the objective 
function is higher than that corresponding to the result provided by the 
first initialization, which is not correct. This fact clearly
indicates the existence of several stationary points or minimizers or both.
Of course it should be nice to establish that the result of the first
initialization is a stationary point and that of the second a minimizer;
but, as already remarked such a verification is practically impossible.
Finally, in the case $m_2 = 17$ also the new initialization is unable to 
provide the correct results.

The results obtained in the multiple image case and described in the next
subsection suggest that this negative result may be due to an insufficient
value of the SNR. Therefore, in the case $m_2 = 17$ we generated an image
which is the sum of 30 frames (we point out that, as already remarked,
in the considered multiple image case we have three times the photons of the
single image case). Using again 2000 iterations, we find that the 
algorithm, with the second initialization, can resolve the binary in the case
SR = 0.81 (even if with a large reconstruction error, about 9 \%, on the PSF) 
but not in the case SR = 0.62. 

However in these difficult cases we observe a new phenomenon: even
if in the limit the results are not satisfactory, the PSF reconstruction error
exhibits a minimum before convergence. If we consider the reconstructions
corresponding to these minima, then, in the case of the first initialization,
the minima do not correspond to a situation where the binary is resolved.
On the other hand, in the case of the second initialization, the binary is
resolved for both SR values, with a 2.03 \% PSF error in the case SR = 0.81
(574 iterations) and a 7.13 \% error in the case SR = 0.62 (1739 iterations).
Such a result presumably indicates the need of introducing a regularization
of the PSF in the objective function, at least for treating the most difficult
cases. In Fig. 3 we show the reconstructions of the PSF 
corresponding to the minimum reconstruction errors. Artifacts due to the 
missed secondary are visible in the case of the first initialization and
also in the case SR = 0.62, since the reconstructed secondary is fainter than
the true one.

\subsubsection{Multiple images}

In this case the integration time of a nonsaturated frame is 95 sec for SR = 
0.77 and 167 sec for SR = 0.46. For each binary and orientation angle we 
consider again 10 frames, so that we have approximately the same number of 
photons in all images. 

We preliminarily remark that, if we compute the value of the normalized 
objective function (which is now given by $2 J_0/3N^2$) by inserting in 
Eq. (\ref{multipleKL}) the noisy and the noise-free images before derotation, 
we expect to obtain a value very close to 1 and this is just what we obtain.
But this is not true if we compute the same quantity using 
the derotated images. Indeed, for the nine binaries as well as for the other
objects, we always obtain a smaller value, namely 0.63. Since this value 
is independent of the object and PSFs, this effect is clearly due to the 
modification of the statistical properties of the data introduced by the 
derotation, as briefly discussed in Sect. 3.2. In any case 
the limiting values of the normalized objective function obtained in our 
experiments are definitely smaller than the values corresponding to the 
input objects and images, an effect already remarked in the previous case.

As in the single image case we first use as initialization a constant array 
for the object and the autocorrelations of the ideal PSFs for the three PSFs.
The results of the reconstructions obtained with this initialization are 
reported in Table 3. 
We obtain that only when both stars have the same magnitude the method is 
able to reconstruct both the binary and the PSFs with sufficient accuracy. 
When we have different magnitudes for the two stars the method is in general 
failing to reproduce the secondary, except in the case of separation 
$d = 80$ mas; in this case a binary with difference of magnitude $\Delta m = 1$ 
is also reconstructed. As in the single image case, the indication 100\% in 
the column for $\Delta m_2/m_2$ means that the method reconstructs an object 
which contains only one bright star (in one case the centroid is shifted one 
pixel in the direction of the secondary.
These results show that, even if we have a higher SNR as already discussed, 
the multiple image case is more difficult than the single one.

\begin{table*}
\caption{Multiple image case - Binary reconstructions provided by the 
algorithm initialized with the autocorrelations of the ideal PSFs. In the 
first column the value of the SR, in the second the angular separation, in 
the third the magnitude of the secondary, in the fourth and fifth the errors 
on the magnitudes of the primary and the secondary star. In the subsequent 
three columns we give the RMSE for the three PSFs. Finally in the last two 
columns we give the value of the normalized objective function, defined by 
$2J_0/3N^2$, as computed at the end of the iterations, and the number of 
outer iterations.}
\label{binaries_m1}
\centering
\begin{tabular}{|c|c|c|cc|ccc|c|c|}
\hline
SR & d (mas) & $m_{2}$ & $\Delta m_1/m_1$ & $\Delta m_2/m_2$ & RMSE$_{0^\circ}$ & RMSE$_{60^\circ}$ & RMSE$_{120^\circ}$ & $J_0$ norm. & IT \\
\hline
  \multirow{9}{*}{0.77} & \multirow{3}{*}{20} & 15 & 0.28\% & 0.23\% & 1.52\% & 2.50\% & 1.87\% & 0.2241 & 1000 \\
  & & 16 & 2.12\% & 100\% & 30.00\% & 31.14\% & 21.79\% & 0.1706  & 1000  \\
  & & 17 & 0.87\% & 100\% & 14.93\% & 15.38\% & 11.00\% & 0.1890  & 1000  \\
\cline{2-10}
  & \multirow{3}{*}{40} & 15 & 0.30\% & 0.28\% & 1.66\% & 1.84\% & 2.56\% & 0.2647  & 1000 \\
  & & 16 & 2.20\% & 100\% & 41.65\% & 33.98\% & 33.80\% & 0.2049  & 1000  \\
  & & 17 & 0.58\% & 100\% & 14.78\% & 14.48\% & 15.90\% & 0.1954  & 1000 \\
\cline{2-10}
  & \multirow{3}{*}{80} & 15 & 0.20\% & 0.21\% & 1.09\%  & 0.83\% & 0.83\% & 0.2180  & 1000  \\
  & & 16 & 0.18\% & 0.23\% & 1.27\% & 0.96\% & 0.99\% & 0.2164  & 1000 \\
  & & 17 & 0.87\% & 100\% & 19.28\% & 19.06\% & 19.04\% & 0.1893  & 1000 \\
\hline
  \multirow{9}{*}{0.46} & \multirow{3}{*}{20} & 15 & 1.27\% & 1.19\% & 9.67\% & 9.69\% & 11.39\% & 0.1335  & 1000  \\
  & & 16 & 0.29\% &  100\% & 32.61\% & 33.10\% & 29.37\% & 0.0836  & 1000 \\
  & & 17 & 0.18\% &  100\% & 13.91\% & 14.12\% & 12.54\% & 0.0795  & 1000 \\
\cline{2-10}
  & \multirow{3}{*}{40} & 15 & 0.89\% & 0.88\% & 5.15\% & 5.62\% & 5.64\% & 0.1516  & 1000 \\
  & & 16 & 0.27\% & 100\% & 47.22\% & 41.54\% & 36.75\% & 0.1360  & 1000 \\
  & & 17 & 0.33\% & 100\% & 13.91\% & 14.67\% & 14.31\% & 0.0944  & 1000 \\
\cline{2-10}
  & \multirow{3}{*}{80} & 15 & 0.68\% & 0.68\% & 3.05\% & 2.60\% & 2.58\% & 0.1042  & 1000 \\
  & & 16 & 0.52\% & 0.60\% & 1.87\% & 1.43\% & 1.44\% & 0.0850   & 1000 \\
  & & 17 & 0.55\% & 100\% & 15.99\% & 15.97\% & 15.98\% & 0.0883  & 1000 \\
  \hline
\end{tabular}
\end{table*}

If we deconvolve the derotated images using the derotated PSFs (this is not 
exactly an {\it inverse crime} because the images were generated with non 
derotated PSFs) all the binaries are correctly reconstructed with small 
errors on the magnitudes. Therefore the failure of our experiment may be 
due again to an inappropriate initialization (the autocorrelations of the ideal 
PSFs have a SR value of about 0.35, much smaller than the SR of the PSFs used 
in image generation) or to inappropriate choices of the internal iterations. 
Also in this case, as in \citet{prato2013} and in the single image case, we
use 50 SGP iterations for the object and one SGP iteration for each PSF.
However several attempts with different numbers of internal iterations did 
not improve the results. Therefore, as in the single image case, we use as
a new initialization of the PSFs the ideal PSFs of LN
with the addition of a small constant selected in such a way to satisfy 
normalization and SR value. The results obtained with this initialization,
using again 50 SGP iterations for the object and one for the PSFs, are
reported in Table 4. Since the convergence is slower than in 
the previous case we use 2000 outer iterations.   

\begin{table*}
\caption{Multiple image case - Binary reconstructions provided by the 
algorithm initialized with the ideal PSFs plus a constant selected 
for satisfying the SR constraint (see the text). The structure of the Table 
is the same of Table \ref{binaries_m1}.}
\label{binaries_m2}
\centering
\begin{tabular}{|c|c|c|cc|ccc|c|c|}
\hline
SR & d (mas) & $m_{2}$ & $\Delta m_1/m_1$ & $\Delta m_2/m_2$ & RMSE$_{0^\circ}$ & RMSE$_{60^\circ}$ & RMSE$_{120^\circ}$ & $J_0$ norm. & IT \\
\hline
  \multirow{9}{*}{0.77} & \multirow{3}{*}{20} & 15 & 0.44\% & 0.34\% & 2.86\% & 4.23\% & 3.42\% & 0.2277 & 2000 \\
  & & 16 & 0.27\% & 0.21\% & 1.56\% & 1.82\% & 1.74\% & 0.2209  & 2000 \\
  & & 17 & 0.07\% & 1.10\% & 2.53\% & 2.70\% & 1.78\% & 0.2095  & 2000 \\
\cline{2-10}
  & \multirow{3}{*}{40} & 15 & 0.45\% & 1.03\% & 5.47\% & 4.61\% & 6.73\% & 0.2670 & 2000\\
  & & 16 & 0.25\% & 0.39\% & 1.63\% & 2.85\% & 2.76\% & 0.2220 & 2000 \\
  & & 17 & 0.11\% & 0.73\% & 2.05\% & 2.78\% & 2.86\% & 0.2102 & 2000 \\
\cline{2-10}
  & \multirow{3}{*}{80} & 15 & 0.35\% & 0.35\% & 2.28\% & 1.51\% & 2.23\% & 0.2204 & 2000 \\
  & & 16 & 0.25\% & 0.26\% & 1.32\% & 1.06\% & 1.14\% & 0.2179 & 2000 \\
  & & 17 & 0.19\% & 0.40\% & 1.32\% & 1.02\% & 0.99\% & 0.2125 & 2000 \\
  \hline
  \multirow{9}{*}{0.46} & \multirow{3}{*}{20} & 15 & 0.80\% & 0.57\% & 4.02\% & 4.51\% & 6.83\% & 0.1037 & 2000 \\
  & & 16 & 0.38\% & 0.95\% & 3.76\% & 3.95\% & 2.52\% & 0.0811 & 2000 \\
  & & 17 & 0.07\% & 6.32\% & 9.05\% & 10.33\% & 6.29\% & 0.0697 & 2000 \\
\cline{2-10}
  & \multirow{3}{*}{40} & 15 & 0.64\% & 2.18\% & 11.23\% & 6.72\% & 8.73\% & 0.1409 & 2000 \\
  & & 16 & 0.49\% & 0.81\% & 2.21\% & 3.20\% & 2.99\% & 0.0837 & 2000 \\
  & & 17 & 0.02\% & 5.89\% & 8.70\% & 7.68\% & 8.84\% & 0.0716 & 2000 \\
\cline{2-10}
  & \multirow{3}{*}{80} & 15 & 0.56\% & 0.55\% & 8.54\% & 4.66\% & 4.64\% & 0.1100 & 2000 \\
  & & 16 & 0.57\% & 0.53\% & 1.99\% & 1.48\% & 1.49\% & 0.0846 & 2000 \\
  & & 17 & 0.48\% & 0.92\% & 2.36\% & 1.95\% & 1.96\% & 0.0785 & 2000 \\
  \hline
\end{tabular}
\end{table*}

\begin{figure*}
  \centering
  \includegraphics[width=0.6\textwidth]{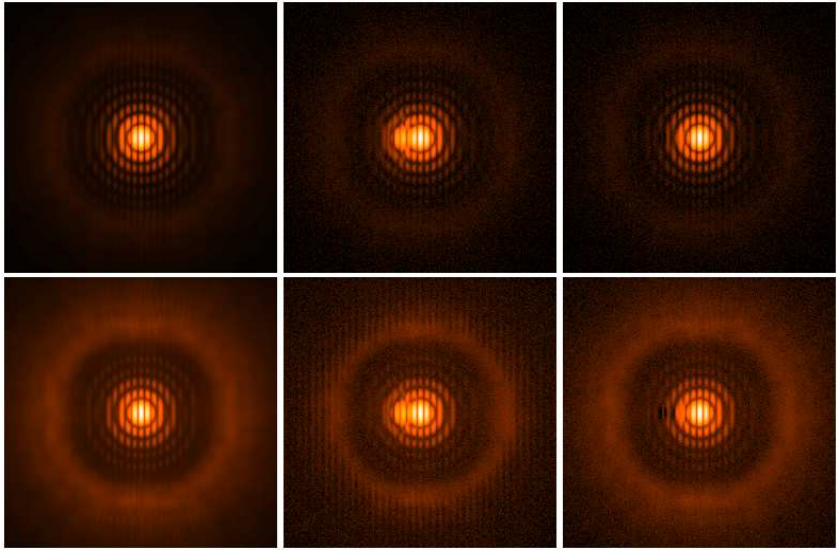} 
  \caption{Multiple image case - PSF reconstruction in the case of the binary 
with $d = 80$ mas and $m_2=17$. First column: the true PSF with SR = 0.77 (top) 
and SR = 0.46 (bottom). Second column: PSF reconstruction provided by the 
algorithm initialized with the autocorrelations of the ideal PSFs. Last column: 
PSF reconstruction if the algorithm is initialized with the ideal PSFs plus a 
constant. All images (only the central part 256$\times$256 is shown) are 
displayed in log scale and correspond to the first orientation of the baseline.}
\label{fig:psf_m}
\end{figure*}

With the new initialization the blind method succeeds in reconstructing all 
the binaries with sufficient accuracy as well as the PSFs. We can add that
in most cases both the normalized objective function and the RMSE on the PSFs 
have a convergent behaviour while, in a few cases, the errors are still 
decreasing after 2000 iterations, thus indicating that a larger number of
iterations could still improve the solution. A comparison of the values of 
the objective function reported in the two tables shows that, in some of 
the cases where the first initialization is failing, the values in Table 
3 are smaller than the corresponding values in Table 
4. This phenomenon was already observed in the single image 
case and means that different stationary points or minimizers are present.

A few more comments on the two tables. If one looks carefully at the reported
results one can remark that, even if the results obtained with the second 
initialization are globally better than those obtained with the first one, 
this may not be true for particular cases (compare, for instance, the results 
for $d = 40$ mas and $\Delta m = 0$). Moreover, the errors obtained with the 
second initialization do not vary in a regular way with the variation of 
angular distance and difference of magnitude. These behaviors can be due
to the fact that 2000 iterations may not be sufficient for assuring
convergence of the method in the case of the second initialization. We
did not push further the iterations because in the case of three 512$\times$512
images the computation time is considerable. By assuming possible fluctuations
due to insufficient number of iterations, a reasonable conclusion seems to
be that, as in the single image case, the two initializations lead to the same 
limit point when the first one is successful.

In Fig. 4 we show an example of reconstructions of the PSF at 
$0^\circ$, for both SR values, when the unknown object is a binary with $d = 
80$ mas and $m_2 = 17$. From the reconstructions displayed in the second column 
and obtained by initializing with the autocorrelations of the ideal PSFs, it 
is evident that they contain a contribution coming from the secondary, while 
this contribution is practically absent in the reconstructions obtained with 
the other initialization and displayed in the third column. 

\begin{figure*}
  \centering
  \includegraphics[width=0.4\textwidth]{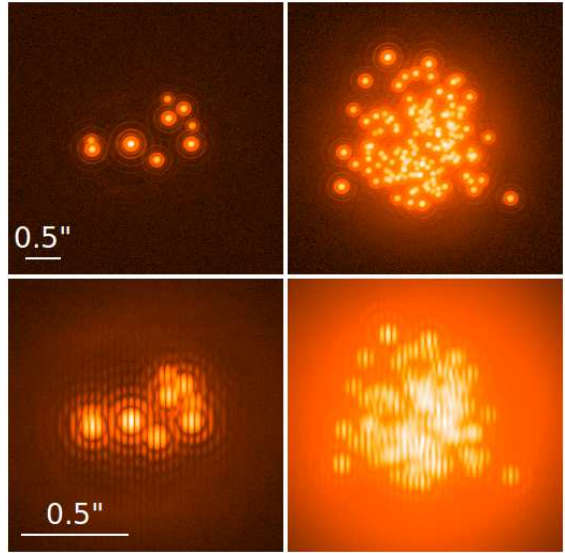} 
  \caption{Top panels: the input images of the ``open star cluster'' (left) and
of the ``globular star cluster'' (right) in the single image case with SR = 0.81. 
Bottom panels: the input images of the two clusters in the case of 
SR = 0.77 and with $0^\circ$ of the baseline (only the central part
256$\times$256 is shown). All images are displayed in log scale. 
The length corresponding to 0.5 arcsec is also indicated.}
  \label{fig:star_clusters}
\end{figure*}

\subsection{Star clusters}

In a second experiment we consider two models of star cluster. The first is
already considered in \citet{prato2013} and is based on an image of the
brightest stars of the Pleiades open cluster; for this reason, we call it
``open star cluster''. It consists of nine stars that we take, in this paper, 
with magnitudes ranging from 14.4 to 17.1. In the single image case, 
the minimum distance between two stars is 120 mas, while the maximum distance 
is 1434 mas, with a mean distance of about 690 mas. In the multiple image 
case, considering the different pixel scale, we reduce of one third all the 
distances.

As a second example we consider a model that we call ``globular star cluster''. 
For simplicity, only 150 stars are considered within the field of view, 
representing a very low crowding condition. The positions of the stars are 
randomly computed following a Gaussian distribution around the center of the 
image (with a standard deviation of about 450 mas in the single image case and 
of about 150 mas in the multiple image case); similarly the magnitudes of the 
stars are randomly distributed around $m$ = 16 with a standard deviation of 
about 0.4. It turns out that the brightest star of the cluster has $m$ = 14.8.

Again, we limit the maximum number of counts in each frame to $5 \times 10^4$, 
keeping fixed to 10 the number of frames.
In Fig. 5 we show the images of the two star clusters
provided by the PSFs with the highest SR. 

\subsubsection{Single image}

\begin{table*}
\caption{Single image case - The reconstruction errors in the case of the two
star cluster models. In the first column the ``open cluster'' is labelled by 
OC, while the ``globular cluster'' is GC. In the second column, we give the 
value of SR, while in the third column we give the initialization of 
the algorithm, denoting by A the autocorrelation of the diffraction-limited 
PSF and by C the diffraction-limited PSF plus a constant selected for 
satisfying the SR constraint (see the text). In the subsequent columns, we give 
the value of the magnitude average reconstruction error 
(MARE) defined in Eq. (\ref{mare}), and the RMSE for the reconstructed PSF. 
Finally in the last two columns we give the value of the normalized objective 
function, defined by $2J_0/N^2$, as computed at the end of the iterations, and 
the number of outer iterations.}
\label{cluster_s}
\centering
\begin{tabular}{|c|c|c|c|c|c|c|}
\hline
 Star Cluster & SR & Init. & MARE & RMSE & $J_0$ norm. & IT \\
\hline
\multirow{4}{*}{OC} & \multirow{2}{*}{0.81} & A & 0.06\% & 0.84\% & 0.5993 & 2000 \\
& & C & 0.06\% & 0.85\% & 0.5991 & 5000 \\
\cline{2-7}
 & \multirow{2}{*}{0.62} & A & 0.09\% & 1.22\% & 0.5165 & 4000 \\
 & & C & 0.09\% & 1.22\% & 0.5165 & 10000 \\
\hline
\multirow{4}{*}{GC} & \multirow{2}{*}{0.81} & A & 0.06\% & 0.87\% & 0.5123 & 3000 \\
& & C & 0.06\% & 0.82\% & 0.4989 & 5000 \\
\cline{2-7}
 & \multirow{2}{*}{0.62} & A & 0.07\% & 1.07\% & 0.4622 & 6000\\
 & & C & 0.07\% & 15.82\% & 0.5698 & 10000 \\
\hline
\end{tabular}
\end{table*}

In the case of the ``open star cluster'', the integration time of a single
frame is 22 sec for SR = 0.81 and 29 sec for SR = 0.62 while in the case of the 
``globular star cluster'' these times are respectively 32 and 42 sec. 

We applied to the four images our blind algorithm using both initializations
introduced in the case of the binaries. The results are reported in Table
5. In the case of the ``open star cluster'' and both values of
SR the two initializations seem to provide sequences of iterations converging 
to the same point. If we look at the image shown in the upper
left panel of Fig. 5 we can observe that it contains
sufficiently well-separated star images which can allow a good estimation
of the PSF by the blind algorithm. 

The situation is a bit different in the case of the ``globular star cluster''
and we can understand this fact if we look at the upper right panel of Fig.
5. In the case of the higher SR value both initializations 
lead essentially to the same result. The small differences may be due to
different convergence rates and could be removed by a more accurate tuning
of the number of iterations. On the other hand in the case of the lower SR
ratio the first initialization, based on the autocorrelation of the 
diffraction-limited PSF, provides the best PSF reconstruction (also
corresponding to a lower value of the objective function). It seems that the
two initializations lead to two different stationary points. In conclusion, 
for this particular object one can state that the first initialization may 
provide a better result than the second one.

\begin{table*}
\caption{Multiple image case - The reconstruction errors in the case of the two 
models of star cluster. The structure is similar to that of Table 
\ref{cluster_s} but now we give the errors on the three PSFs and the
normalized objective function is defined by $2J_0/3N^2$.}
\label{cluster_m}
\centering
\begin{tabular}{|c|c|c|c|ccc|c|c|}
\hline
 Star Cluster & SR & Init. & MARE & RMSE$_{0^\circ}$ & RMSE$_{60^\circ}$ & RMSE$_{120^\circ}$ & $J_0$ norm. & IT \\
\hline
\multirow{4}{*}{OC} & \multirow{2}{*}{0.77} & A & 0.35\% & 2.14\% & 4.28\% & 4.19\% & 0.3049 & 1000 \\
& & C & 0.59\% & 4.16\% & 7.67\% & 7.62\% & 0.2997 & 5000 \\
\cline{2-9}
 & \multirow{2}{*}{0.46} & A & 0.81\% & 3.43\% & 6.07\% & 6.00\% & 0.1321 & 2000 \\
 & & C & 0.89\% & 3.74\% & 7.77\% & 7.81\% & 0.1237 & 10000 \\
\hline
\multirow{4}{*}{GC} & \multirow{2}{*}{0.77} & A & 0.38\% & 1.98\% & 3.46\% & 3.38\% & 0.7597 & 3000 \\
& & C & 0.71\% & 11.00\% & 11.06\% & 12.97\% & 0.7459 & 5000 \\
\cline{2-9}
 & \multirow{2}{*}{0.46} & A & 1.04\% & 5.63\% & 9.69\% & 9.38\% & 0.3557 & 6000\\
 & & C & 0.90\% & 25.81\% & 16.37\% & 17.94\% & 0.3043 & 10000 \\
\hline
\end{tabular}
\end{table*}

\subsubsection{Multiple images}

In the case of the ``open cluster'' model, the integration time is 53 sec for 
SR = 0.77 and 93 sec for SR = 0.46. On the other hand the integration time for 
the ``globular cluster'' images is 78 sec for SR = 0.77 and 136.5 sec for 
SR = 0.46. 

In both cases we apply our blind algorithm using the two initializations 
already used in the previous sections, with 50 inner SGP iterations for 
the object and one iteration for each PSF. The results obtained for the
``open cluster'' with the two initializations are given in the 
first two rows of Table 6 in the case SR = 0.77 and in the
following two rows those obtained in the case SR = 0.46.
Similarly the results obtained for the ``globular cluster'' are given in 
the second half of the same table. 

In the multiple image case the situation is more complex than in the single
one, and this is not surprising since now we must reconstruct four blocks
of variables. By looking at the results reported in Table 6
we see that the two initializations produce in all cases two sequences of
iterations converging to distinct results. Even if, in some cases, the two 
values of the objective function are very close, the corresponding points
are definitely different, thus implying the existence of several minimizers
or stationary points with very close values of the objective functions.

It is interesting to remark that, while in the case of the binaries the best
results are provided by the second initialization, now they are provided
by the first one, based on the autocorrelations of the ideal PSFs. The highest 
reconstruction errors are obtained in the case of the lowest SR, as one 
should expect.
We also remark that in the case of the second initialization we used a larger
number of iterations because the convergence is slower than in the case of
the first initialization. From the comparison of the results obtained for 
the binaries with those obtained for the star clusters we deduce that the
problem of the initial PSFs is essentially open; therefore, in the case of
practical applications, one should try with different initializations,
using also physical intuition in their choice.

As a final comment, all the values of the objective function corresponding to
the best solutions are higher than those corresponding to the other ones.

\section{Conclusions}

In this paper we extend to the case of Fizeau interferometry a blind
deconvolution method previously proposed for single aperture telescopes and we
validate the method in both cases, called respectively multiple image
and single image case. 

It is well-known that the problem of blind deconvolution is extremely ill-posed
and the introduction of constraints on PSF and object does not exclude the
existence of several local minima, stationary points etc. In our approach
the most significant constraint is the SR constraint on the PSFs, as 
suggested in \citet{desidera2009}. This constraint excludes the trivial
solution of a delta function for the PSF and image for the object.

From our numerical analysis it turns out that the problem of Fizeau
interferometry is more difficult than the problem of single aperture
telescopes. The reason may be twofold. On one hand the number of variables
to be reconstructed is larger and it is known that in the minimization
of a nonconvex block-convex function the theoretical results are weaker
when the number of blocks is greater than two. On the other hand the
PSFs are very structured due to the presence of interference fringes so
that if the initialization does not contain sufficient information on these 
structures it is difficult if not impossible to obtain acceptable results.

An astonishing feature already observed in the single aperture case and
confirmed in the present paper is that very often the value of the objective
function corresponding to a sensible solution is greater than the value 
corresponding to an unacceptable one. Obviously it is impossible to
verify if these points correspond to local minima or to stationary points.
In any case this result raises the issue if global minima, in case
they were computable, provide sensible solutions or not.

In summary the results of this paper open a large number of problems;
however we think that the proposed method, which has a sound mathematical
foundation, is very flexible and can help to investigate these problems. 
Moreover, last but not least, it can be
extended to introduce regularization terms both on the object and on the
PSF (or PSFs) thanks to the high flexibility of SGP, which is the basic tool 
in our approach. It is sufficient to modify the scaling
factor along the lines suggested in \citet{lanteri}. Obviously, in such a
case, the additional problem arises of the choice of the regularization
parameters.

We conclude by remarking that the IDL routines implementing our method are 
available on request and will be included in the forthcoming version of the
Software Package AIRY (v.6.1) \citep{carbillet2014}.

\section*{Acknowledgements}

This work has been partially supported by MIUR (Italian Ministry for University 
and Research), under the projects FIRB - Futuro in Ricerca 2012 (contract 
RBFR12M3AC) and PRIN 2012 (contract 2012MTE38N), and by INAF (National Institute for 
Astrophysics) under the project TECNO-INAF 2010 ``Exploiting the adaptive 
power: a dedicated free software to optimize and maximize the scientific output
of images from present and future adaptive optics facilities''. The Italian 
GNCS - INdAM (Gruppo Nazionale per il Calcolo Scientifico - Istituto Nazionale 
di Alta Matematica) is also acknowledged.

\section*{References}

\bibliographystyle{elsarticle-harv}
\bibliography{paper}

\end{document}